\documentclass{amsart}
\usepackage{a4wide,amsfonts,amssymb,stmaryrd,amscd,amsmath,latexsym,amsbsy}
\numberwithin{equation}{section}
\theoremstyle{plain}
\newtheorem{thm}{Theorem}[section]

\newtheorem{lem}[thm]{Lemma}

\theoremstyle{remark}
\newtheorem{rema}[thm]{Remark}
\title{Generalized Cherednik-Macdonald identities}
\author{Jasper V. Stokman}
\address{Korteweg-de Vries Institute for Mathematics,
University of Amsterdam, Plantage Muidergracht 24, 1018 TV Amsterdam,
The Netherlands}
\email{jstokman@science.uva.nl}
\begin{document}
\begin{abstract}
We derive generalizations of the Cherednik-Macdonald constant
term identities associated to root systems which depend, besides
on the usual multiplicity function, symmetrically on two
quasi-periods $\omega_{\pm}$. They are natural analogues of the Cherednik-Macdonald constant term 
$q$-identities in which 
the deformation parameter $q=\exp(2\pi i\omega_+/\omega_-)$ 
is allowed to have modulus one.
They unite the 
Cherednik-Macdonald constant term $q$-identities 
with closely related Jackson $\widetilde{q}$-integral identities 
due to Macdonald, where the deformation parameter 
$\widetilde{q}=\exp(-2\pi i\omega_-/\omega_+)$ is related to $q$ by modular inversion.
\end{abstract}

\maketitle


\section{Introduction}

The Cherednik-Macdonald identities, conjectured by Macdonald \cite{MacC} and proven in full generality by Cherednik \cite{C}, are explicit constant term evaluations for
certain families of densities $\Delta$ on compact tori. The density $\Delta$ (see \eqref{Delta})
depends on a root system $\Sigma$ in Euclidean space $V$, on a multiplicity label $k$ on $\Sigma$
(free parameters) and on a deformation parameter $q=\exp(2\pi i\omega_+/\omega_-)$ satisfying $|q|<1$.
It is viewed as a density on the compact torus $T=V/P^\vee$, where $P^\vee$ is the co-weight lattice of $\Sigma$ in $V$, and can
be explicitly expressed in terms of $q$-gamma functions. 
In his well known 1987 manuscript Macdonald \cite{Mac1} showed that $\Delta$, viewed as density on the slightly enlarged torus
$V/Q^\vee$ with $Q^\vee$ the co-root lattice of $\Sigma$ in $V$, serves
as the weight function for a remarkable set of multivariate orthogonal polynomials, nowadays known as the Macdonald polynomials.
Cherednik's \cite{C} proof of the constant term identities highlights the crucial role of the double affine Hecke algebra in the theory of Macdonald polynomials.

In this paper we derive natural analogues of the Cherednik-Macdonald
constant term identities for $|q|=1$. They are expected to naturally appear in harmonic analysis on certain non-compact quantum groups
(see e.g. \cite {PT, B, R1, PW} for compelling evidence in this direction) and in certain classes of integrable systems
(compare e.g. with \cite{JM, BT, FKV, KLS}). 

The appearance of the $q$-gamma function in the Cherednik-Macdonald identities 
is the first apparent obstacle for the generalization of the identities to $|q|=1$, since the $q$-gamma function is only well defined for $|q|<1$.
The key to overcome this hurdle is Shintani's \cite{Sh} observation
that a suitable quotient of the $q$-gamma function and the
$\widetilde{q}$-gamma function, where $\widetilde{q}=\exp(-2\pi i\omega_-/\omega_+)$ 
is the deformation parameter related to $q=\exp(2\pi i\omega_+/\omega_-)$ by modular inversion, 
admits analytic continuation to a parameter regime containing co-linear quasi-periods $\omega_+$ and $\omega_-$ (i.e. $|q|=|\widetilde{q}|=1$).
The analytic continuation of this remarkable quotient serves as the natural analogue of the $q$-gamma function for $|q|=1$. It has (re)appeared in various guises and in
different contexts over the past century and goes by names as the double gamma function (Barnes \cite{Ba}),
the double sine function (Kurokawa \cite{Ku}), the quantum exponential function (Woronowicz \cite{W}), the hyperbolic gamma function (Ruijsenaars \cite{R0}), 
the noncompact quantum dilogarithm (Faddeev and Kashaev, see e.g. \cite{FKV}) and the $\gamma$-function (Volkov \cite{V}). 
In this paper we have chosen to relate it to Ruijsenaars' \cite{R0} hyperbolic gamma function (through the formulas \eqref{qshift} and \eqref{ShG} below).
The precise links with the other close relatives can be easily deduced from the appendices in \cite{R, PT}. 

We thus seek a generalization of the Cherednik-Macdonald constant term identities in which the role of the $q$-gamma function is taken over by
the hyperbolic gamma function. To achieve this, we unite the 
constant terms with absolutely convergent sums of the form
\[\sum_{\lambda\in P^\vee}\widetilde{\Delta}(v+\lambda),
\] 
for suitable dual densities $\widetilde{\Delta}$ (see \eqref{Delta}) that can be explicitly expressed in terms of $\widetilde{q}$-gamma functions
(alternatively, the sums can be written as multidimensional Jackson $\widetilde{q}$-integrals).
The sums are independent of $v$ and can be explicit evaluated using Macdonald's \cite{Mac2} summation identities. 
These identities are in some sense "dual" summation versions of the Cherednik-Macdonald constant term identities.
Combining both of them gives an explicit evaluation
of the absolutely convergent integral
\begin{equation}\label{introintegral}
\int_V\Delta(v)\widetilde{\Delta}(v)dv=\int_D\Delta(v)\Bigl(\sum_{\lambda\in P^\vee}\widetilde{\Delta}(v+\lambda)\Bigr)dv
\end{equation}
with $dv$ Lebesgue measure on $V$ and 
$D\subset V$ a fundamental domain for the translation action of $P^\vee$ on $V$. 
It turns out that the density $\Delta\widetilde{\Delta}$ (as well as the explicit evaluation of the integral \eqref{introintegral})
can be entirely expressed in terms of hyperbolic gamma functions. Furthermore,
after a suitable rotation of the integration region $V$ within its complexification $V_{\mathbb{C}}$, 
the resulting integral identity admits analytic continuation to a parameter regime in which the two quasi-periods $\omega_+$ and $\omega_-$
are allowed to be co-linear (in which case $|q|=|\widetilde{q}|=1$). The identities thus obtained are the generalized Cherednik-Macdonald identities
referred to in the title of the paper.

Integral identities involving hyperbolic gamma functions (for which we will use the terminology hyperbolic integral identities), have 
appeared at various places in the literature. There is a large supply of univariate hyperbolic integral evaluations, see e.g. \cite{PT, R1, S, V, BRS}. In particular, 
the method employed in this paper was used in \cite{S} to obtain univariate hyperbolic beta integrals (containing as special case
the generalized Cherednik-Macdonald identity for $\Sigma$ of rank one). 
Multivariate hyperbolic integral evaluations have been obtained in e.g. \cite{Rains, vDS, Sp, BR}. The multivariate hyperbolic integral evaluation relevant to
the generalized Cherednik-Macdonald identities is the type II multivariate hyperbolic integral evaluation from \cite[Thm. 4]{vDS} (see also \cite[Cor. 4.4]{Rains}). 
Concretely, its first level degeneration 
(see \cite[Thm. 5]{vDS}) 
is the generalized constant term identity associated to Koornwinder's \cite{K} extension of the Macdonald theory 
(in this case the associated root system $\Sigma$ is the nonreduced root system of type $BC$). It contains the generalized Cherednik-Macdonald
identities associated to root systems of type $A_1$, $B_n$ and $C_n$ as special cases.
The techniques of the present paper can in fact easily be generalized to include the Koornwinder case, but we do
not pursue this here in detail.

The contents of the paper is as follows.
After a short introduction on the hyperbolic gamma function in subsection \ref{S21} 
and fixing notations on root systems in subsection \ref{S22}
we formulate the generalized Cherednik-Macdonald identities in subsection \ref{S23}. In subsection \ref{S31}
we fix $\omega_+/\omega_-$ in the upper half plane (so that $|q|, |\widetilde{q}|<1$), we use 
the expression of the hyperbolic gamma function as quotient of $q$-gamma and $\widetilde{q}$-gamma functions
and rotate the integration regime in order to be able to fold the integral as in \eqref{introintegral}. 
In subsection \ref{S32} we recall the Cherednik-Macdonald constant term identities and the associated dual
summation identities to complete the proof of the generalized Cherednik-Macdonald identities. In section 4
we show that the generalized Cherednik-Macdonald identities associated to root systems $\Sigma$ of type $A_1$, $B_n$
and $C_n$ are special cases of the first level degeneration \cite[Thm. 5]{vDS} of the multivariate hyperbolic
integral evaluation of type II.\\

\noindent
{\it Convention:} We take the branch of $\sqrt{\cdot}$ which is nonnegative on
the nonnegative real axis and with branch cut along the negative real axis.\\

\noindent
{\bf Acknowledgments:} The author is supported by the Netherlands Organization for Scientific
Research (NWO) in the VIDI-project "Symmetry and modularity in exactly solvable models".
He thanks Fokko van de Bult and Eric Rains for stimulating discussions.


\section{Formulation of the integral identities}\label{Section2}


\subsection{The hyperbolic gamma function}\label{S21}

We consider Cherednik-Macdonald constant term identities in which the role of the $q$-shifted factorial
\begin{equation}\label{qshift}
\bigl(z;q\bigr)_{\infty}:=\prod_{j=0}^{\infty}\bigl(1-q^jz\bigr),\qquad |q|<1
\end{equation}
is replaced by Ruijsenaars' \cite{R0} hyperbolic gamma function (the $q$-gamma function $\Gamma_q$ referred to in the introduction
relates to the $q$-shifted factorial by $\Gamma_q(x)=(1-q)^{1-x}\bigl(q;q\bigr)_{\infty}/\bigl(q^x;q\bigr)_{\infty}$ for $0<q<1$).
All results stated here can be immediately traced back to \cite{R0} and \cite[Appendix A]{R}. See also the paper \cite{V}, which gives a nice
overview of some of the key properties of the closely related $\gamma$-function.

We write
\[\mathbb{C}_{\pm}=\{z\in\mathbb{C} \,\, | \,\, \hbox{Re}(z)\gtrless 0\}
\]
for the open right/left half-plane in $\mathbb{C}$ and
\[
\mathbb{H}_{\pm}=\{z\in\mathbb{C} \,\, | \,\, \hbox{Im}(z)\lessgtr 0\}
\]
for the open upper/lower half-plane in $\mathbb{C}$.
The hyperbolic gamma function $G(\omega_+,\omega_-;z)$, depending on two quasi-periods $\omega_{\pm}\in\mathbb{C}_+$, is defined
for $|\hbox{Im}(z)|<\frac{1}{2}\hbox{Re}(\omega_++\omega_-)$ by
\[G(\omega_+,\omega_-;z)=\exp\left(i\int_{0}^\infty\frac{dy}{y}\left(\frac{\sin(2yz)}{2\sinh(\omega_+y)\sinh(\omega_-y)}-
\frac{z}{\omega_+\omega_-y}\right)\right).
\]
We suppress the quasi-periods and write $G(z)=G(\omega_+,\omega_-;z)$ if no confusion can arise. We denote $\omega=\frac{1}{2}(\omega_++\omega_-)$
and $\Lambda=\Lambda_{\omega_+,\omega_-}:=\mathbb{Z}_{\geq 0}i\omega_++\mathbb{Z}_{\geq 0}i\omega_-\subset\overline{\mathbb{H}}_+$ (where $\overline{\mathbb{H}}_+$ is the
closed upper half plane in $\mathbb{C}$). The hyperbolic gamma function satisfies the functional equations
\begin{equation}\label{functionalequation}
\frac{G(z+i\omega_{\pm}/2)}{G(z-i\omega_{\pm}/2)}=2\cosh(\pi z/\omega_{\mp})
\end{equation}
whenever the left hand side is defined. The functional equations allows to extend $G(\omega_+,\omega_-;z)$ to a meromorphic function on 
$(\omega_+,\omega_-,z)\in\mathbb{C}_+^{\times 2}\times\mathbb{C}$ with zero and singular locus contained in
$\{(\omega_+,\omega_-,z)\in\mathbb{C}_+^{\times 2}\times\mathbb{C} \, | \, z\in i\omega+\Lambda_{\omega_+,\omega_-} \}$ 
and $\{(\omega_+,\omega_-,z)\in\mathbb{C}_+^{\times 2}\times\mathbb{C} \, | \, z\in -i\omega-\Lambda_{\omega_+,\omega_-}\}$, respectively. 

Important for our considerations is the fact that $G(z)$ may be viewed as a generalization of the $q$-shifted factorial 
\eqref{qshift} in which $q$ is allowed to have modulus one. This
becomes transparent from Shintani's \cite{Sh} product formula
\begin{equation}\label{ShG}
G(z)=\frac{\bigl(\exp(-2\pi(z-i\omega)/\omega_-);q\bigr)_{\infty}}{\bigl(\exp(-2\pi(z+i\omega)/\omega_+);\widetilde{q}\bigr)_{\infty}}
\exp\left(-\frac{\pi i}{24}\Bigl(\frac{\omega_+}{\omega_-}+\frac{\omega_-}{\omega_+}\Bigr)\right)\exp\left(-\frac{\pi iz^2}{2\omega_+\omega_-}\right)
\end{equation}
for $\omega_+/\omega_-\in\mathbb{H}_+$, where the bases $q$ and $\widetilde{q}$ are related to the quasi-periods
$\omega_{\pm}$ by
\[q=\exp(2\pi i\omega_+/\omega_-),\qquad \widetilde{q}=\exp(-2\pi i\omega_-/\omega_+).
\]
For a simple proof of \eqref{ShG}, see e.g. \cite[Prop. 6.1]{S}.
Note that the requirement $\omega_+/\omega_-\in\mathbb{H}_+$ is necessary for the
right hand side of \eqref{ShG} to be well defined, since it 
implies that $|q|<1$ and $|\widetilde{q}|<1$.
On the other hand, $|q|=1$ and $|\widetilde{q}|=1$
corresponds to co-linear quasi-periods $\omega_{\pm}\in\mathbb{C}_+$, in which case the hyperbolic gamma function $G(z)$ itself still makes
perfect sense. 


\subsection{Root systems}\label{S22}
The Cherednik-Macdonald constant terms and their generalizations in the present paper are naturally attached to root systems. We fix in this
subsection the necessary notations. Let $\Sigma\subset V$ be an irreducible, crystallographic, reduced root system in an Euclidean space
$V$ of rank $n$. We write $\langle\cdot,\cdot\rangle$ for the scalar product of $V$ and $\|\cdot\|$
for the associated norm on $V$. We normalize the root system $\Sigma$ such that $\|\alpha\|^2=2$
for short roots $\alpha\in\Sigma$. The co-roots are $\alpha^\vee=2\alpha/\|\alpha\|^2$ for $\alpha\in\Sigma$.
Let $W$ be the Weyl group of $\Sigma$, and denote $Q$, $Q^\vee$, $P$ and $P^\vee$ for the root lattice, co-root lattice, weight lattice
and co-weight lattice of $\Sigma$, respectively.
Note that $Q$ (respectively $Q^\vee$) is a sublattice of $P$ (respectively $P^\vee$)
of finite index. We write
\begin{equation}\label{f}
f=\#(P/Q)
\end{equation}
for the index of $Q$ in $P$. It is also equal to the index of $Q^\vee$ in $P^\vee$.

We fix a basis $\{\alpha_j\}_{j=1}^n$ for $\Sigma$. We denote
$\Sigma^+$ and $\Sigma^-=-\Sigma^+$ for the associated sets of
positive and negative roots in $\Sigma$, and
we write $\omega_j\in P$ and $\omega_j^\vee\in P^\vee$ ($1\leq j\leq n$)
for the associated fundamental weights and co-weights.

A $W$-invariant complex valued function $k: \Sigma\rightarrow \mathbb{C}$ is called
a multiplicity function. We write $k_\alpha$ for the value of
the multiplicity function $k$ at $\alpha\in\Sigma$. Let $\mathcal{K}$ 
be the complex vector space consisting of multiplicity functions. It is one dimensional if
all roots $\alpha\in\Sigma$ have the same length and two dimensional otherwise.
For $k\in\mathcal{K}$ we write $\rho_k=\frac{1}{2}\sum_{\alpha\in \Sigma^+}k_\alpha\alpha$
for the $k$-deformation of the half sum of positive roots. Note that
\begin{equation}\label{rhofundamental}
\rho_k=\sum_{j=1}^nk_{\alpha_j}\omega_j.
\end{equation}

Following Macdonald \cite{Mac1}, \cite{Mac2}, we formulate two types of generalized Cherednik-Macdonald
integral identities associated to $\Sigma$. The two types depend
on two different choices {\bf (i)} and {\bf (ii)} of an auxiliary multiplicity function $u\in\mathcal{K}$:
\begin{enumerate}
\item[{{\bf (i)}}] $u_\alpha=1$ for all $\alpha\in\Sigma$.
\item[{{\bf (ii)}}] $u_\alpha=2/\|\alpha\|^2$
for all $\alpha\in\Sigma$.
\end{enumerate}
The possible values of $u_\alpha$ are $1,\frac{1}{2}$ or $\frac{1}{3}$. We write $\alpha^\prime:=u_\alpha\alpha$ for $\alpha\in\Sigma$,
so that $\alpha^\prime=\alpha$ for case {\bf (i)} and $\alpha^\prime=\alpha^\vee$ for case {\bf (ii)}.

\subsection{The integral identities}\label{S23}

For $\alpha\in\Sigma$ and quasi-periods $\omega_\pm\in\mathbb{C}_+$ we write
\[G_\alpha(z)=G_\alpha(\omega_+,\omega_-;z):=G(\omega_+,u_\alpha\omega_-;z),
\]
$\omega_\alpha=\frac{1}{2}(\omega_++u_\alpha\omega_-)$ and $\Lambda_\alpha=\mathbb{Z}_{\geq 0}i\omega_++\mathbb{Z}_{\geq 0}iu_\alpha\omega_-$.
Since $G(r\omega_+,r\omega_-;rz)=G(\omega_+,\omega_-;z)$ for $r>0$ and $G(\omega_+,\omega_-;z)=G(\omega_-,\omega_+;z)$ (both are immediate from the
definition of $G(z)$) 
we can alternatively write $G_\alpha(z)=G(\omega_-, u_\alpha^{-1}\omega_+;u_\alpha^{-1}z)$.

The integral identities depend, besides on the root system $\Sigma$, on two quasi-periods $\omega_{\pm}$ and on a multiplicity function $k:\Sigma\rightarrow\mathbb{C}$.
We formulate the integral identities for parameters in the open, arcwise connected parameter space
\[\mathcal{S}=\{(\omega_+,\omega_-,k)\in\mathbb{C}_+^{\times 2}\times\mathcal{K} \,\, | \,\,
\omega_+\omega_-\in\mathbb{H}_-,\quad k_\alpha\in\mathbb{C}_-\cap\omega_+\omega_-\mathbb{C}_+\,\quad \forall\,\alpha\in\Sigma\}.
\]
So for quasi-periods $\omega_{\pm}\in\mathbb{C}_+$ satisfying $\omega_+\omega_-\in\mathbb{H}_-$ 
we have $(\omega_+,\omega_-,k)\in\mathcal{S}$ for $k\in\mathcal{K}$ iff
$k$ takes value in the non-empty open wedge within $\mathbb{C}_-$
bounded by the half-lines $\omega_+\omega_-i\mathbb{R}_{\leq 0}$ and $i\mathbb{R}_{\leq 0}$.

Let $dv$ be Lebesgue measure on $V$ normalized by $\int_Ddv=1$, where $D$ is the parallelepiped in $V$ spanned by the fundamental co-weights
$\omega_j^\vee$ ($1\leq j\leq n$).
The extension of the scalar product $\langle \cdot, \cdot\rangle$ to a complex bilinear form on the complexification $V_{\mathbb{C}}$
of $V$ will also be denoted by $\langle \cdot,\cdot\rangle$.
\begin{thm}\label{hyperbolicCM}
For $(\omega_+,\omega_-,k)\in\mathcal{S}$ we have
\begin{equation*}
\begin{split}
\int_V\prod_{\alpha\in\Sigma}&\frac{G_\alpha(\langle\alpha^\prime,v\rangle+i\omega_\alpha)}
{G_\alpha(\langle\alpha^\prime,v\rangle+i(k_\alpha+\omega_\alpha))}dv=\\
&\quad=f\#W\prod_{j=1}^n\sqrt{\frac{\omega_+\omega_-}{u_{\alpha_j}}}
\prod^\prime_{\alpha\in\Sigma^+}{}\frac{G_\alpha(i(\langle\rho_k,\alpha^\vee\rangle+\omega_\alpha))G_\alpha(i(\langle\rho_k,\alpha^\vee\rangle-\omega_\alpha))}
{G_{\alpha}(i(\langle\rho_k,\alpha^\vee\rangle+k_\alpha+\omega_\alpha))G_\alpha(i(\langle\rho_k,\alpha^\vee\rangle-k_\alpha-\omega_\alpha))},
\end{split}
\end{equation*}
where the regularized product $\overset{\prime}{\prod}$ means that the 
factors $G_{\alpha_j}(i(\langle\rho_k,\alpha_j^\vee\rangle-k_{\alpha_j}-\omega_{\alpha_j}))$ \textup{(}$1\leq j\leq n$\textup{)}
in the denominator should be omitted.
\end{thm}

\begin{rema}
Using the $W$-invariance of the integrand, \eqref{rhofundamental}, and the special value (cf. \cite[A.8]{S})
\[
G\Bigl(\omega_+,\omega_-;\frac{i}{2}(-\omega_++\omega_-)\Bigr)=\sqrt{\omega_-/\omega_+}
\]
of the hyperbolic gamma function, the integral identity is equivalent to
\begin{equation}\label{CMalternative}
\begin{split}
\int_{V^+}\prod_{\alpha\in\Sigma}&\frac{G_\alpha(\langle\alpha^\prime,v\rangle+i\omega_\alpha)}
{G_\alpha(\langle\alpha^\prime,v\rangle+i(k_\alpha+\omega_\alpha))}dv=\\
&=f\omega_-^n\prod_{\alpha\in\Sigma^+}\frac{G_\alpha(i(\langle\rho_k,\alpha^\vee\rangle+\omega_\alpha))G_\alpha(i(\langle\rho_k,\alpha^\vee\rangle-\omega_\alpha))}
{G_{\alpha}(i(\langle\rho_k,\alpha^\vee\rangle+k_\alpha+\omega_\alpha))
G_\alpha(i(\langle\rho_k,\alpha^\vee\rangle-k_\alpha+u_\alpha\omega_-\delta_\alpha-\omega_\alpha))}
\end{split}
\end{equation}
where 
\[V^+:=\{v\in V \,\, | \,\, \langle\alpha,v\rangle\geq 0 \quad \forall\,\alpha\in\Sigma^+ \}=\bigoplus_{j=1}^n\mathbb{R}_{\geq 0}\omega_j^\vee
\]
is the closed positive Weyl chamber and 
$\delta_\alpha$ ($\alpha\in\Sigma^+$) is $1$ if $\alpha\in\{\alpha_1,\ldots,\alpha_n\}$
and is $0$ otherwise. See also section \ref{s4} for alternative expressions of the integral identities for root
systems $\Sigma$ of type $A_1$, $B_n$ and $C_n$.
\end{rema}
We now verify that both sides of \eqref{CMalternative} are well defined, and that they depend
analytically on $(\omega_+,\omega_-,k)\in\mathcal{S}$. The actual proof of the integral identity \eqref{CMalternative} is postponed
to the next section.
Write
\begin{equation}\label{integrand}
I(v):=\prod_{\alpha\in\Sigma}\frac{G_\alpha(\langle\alpha^\prime,v\rangle+i\omega_\alpha)}
{G_\alpha(\langle\alpha^\prime,v\rangle+i(k_\alpha+\omega_\alpha))}
\end{equation}
for the integrand of \eqref{CMalternative}, viewed as meromorphic function on $v\in V_{\mathbb{C}}$.

\begin{lem}
{\bf (a)} For parameters $(\omega_+,\omega_-,k)\in\mathcal{S}$ the integral $\int_{V^+}I(v)dv$ is absolutely convergent. It depends
analytically on $(\omega_+,\omega_-,k)\in\mathcal{S}$.\\
{\bf (b)} The right hand side of \eqref{CMalternative} depends analytically on $(\omega_+,\omega_-,k)\in\mathcal{S}$.
\end{lem}
\begin{proof}
{\bf (a)} First we show that the integrand $I(v)$ is analytic at $v\in V$. For this it is convenient to rewrite the integrand $I(v)$ using
the reflection equation \cite[Prop. III.2]{R0}
\begin{equation}\label{reflectionequation}
G(z)G(-z)=1
\end{equation}
for the hyperbolic gamma function (its validity is immediate from the definition of $G(z)$). 
Together with \eqref{functionalequation} it yields  
\begin{equation}\label{analyticproduct}
G(z+i\omega)G(-z+i\omega)=4\sinh(\pi z/\omega_+)\sinh(\pi z/\omega_-).
\end{equation}
Using \eqref{reflectionequation} and \eqref{analyticproduct} the integrand $I(v)$
can then be rewritten as 
\begin{equation}\label{integrandalternative}
\begin{split}
I(v)=\prod_{\alpha\in\Sigma^+}&\left\{4\sinh\bigl(\pi\langle\alpha^\prime,v\rangle/\omega_+\bigr)
\sinh\bigl(\pi\langle\alpha,v\rangle/\omega_-\bigr)\right.\\
&\qquad\left.\times G_\alpha\bigl(\langle\alpha^\prime,v\rangle-i(k_\alpha+\omega_\alpha)\bigr)
G_\alpha\bigl(-\langle\alpha^\prime,v\rangle-i(k_\alpha+\omega_\alpha)\bigr)\right\},
\end{split}
\end{equation}
from which it immediately follows that the possible poles of $v\mapsto I(v)$ are at
\begin{equation}\label{polesI}
\langle\alpha^\prime,v\rangle\in -ik_\alpha+\Lambda_\alpha,\qquad \alpha\in\Sigma.
\end{equation}
Since $(\omega_+,\omega_-,k)\in\mathcal{S}$ we have $k_\alpha\in\mathbb{C}_-$ and $\Lambda_\alpha\subset \overline{\mathbb{H}}_+$,
hence $-ik_\alpha+\Lambda_\alpha\subset\mathbb{H}_+$  for all $\alpha\in\Sigma$. It follows that $I(v)$ is analytic at $v\in V$. 

For the convergence of the integral \eqref{CMalternative} we use asymptotic estimates for the hyperbolic gamma function from \cite[Thm. A.1]{R} 
(see also \cite{Rains}), which imply that for compacta $K_{\pm}\subset\mathbb{C}_+$, $K\subset\mathbb{R}$ there exist
$R,C>0$ depending only on $K_{\pm}$ and $K$ such that
 \[
|G(\omega_+,\omega_-;z)|\leq C|\exp\bigl(\mp\pi iz^2/2\omega_+\omega_-\bigr)|
\]
when $\hbox{Re}(z)\gtrless R$, $\hbox{Im}(z)\in K$ and $\omega_{\pm}\in K_{\pm}$. Applied to the alternative expression
\eqref{integrandalternative} of the integrand $I(v)$, we obtain for compacta $K\subset \mathcal{S}$ the estimate
\[|I(v)|\leq C_K|\exp\bigl(-4\pi\langle\rho_k,v\rangle/\omega_+\omega_-\bigr)|,\qquad \forall\, v\in V^+,\quad
\forall\, (\omega_+,\omega_-,k)\in K
\]
for some constant $C_K>0$. Observing that $-k_\alpha/\omega_+\omega_-\in\mathbb{C}_-$ ($\alpha\in\Sigma$) if $(\omega_+,\omega_-,k)\in\mathcal{S}$
and that $2\langle\rho_k,\omega_j^\vee\rangle$ is a non-empty sum of $k_\alpha$'s for $1\leq j\leq n$, 
it now easily follows that $\int_{V^+}I(v)dv$ is absolutely convergent and that it depends analytically on
$(\omega_+,\omega_-,k)\in \mathcal{S}$.\\
{\bf (b)}
For $\alpha\in\Sigma^+\setminus \{\alpha_1,\ldots,\alpha_n\}$,
\[\langle\rho_k,\alpha^\vee\rangle-k_\alpha\]
is a nonempty sum of $k_\beta$'s. Furthermore, for $(\omega_+,\omega_-,k)\in\mathcal{S}$ we have $\Lambda_\alpha\subset \overline{\mathbb{H}}_+\cap\omega_+\omega_-\overline{\mathbb{H}}_+$.
Using these two observations it is straightforward to check that the right hand side of \eqref{CMalternative} depends analytically on $(\omega_+,\omega_-,k)\in\mathcal{S}$. 
\end{proof}
\begin{rema}\label{T}
The proof of the lemma shows that the conditions $k_\alpha\in\mathbb{C}_-$ are needed for the singular locus of the integrand $I(v)$
to be properly separated by the integration region $V$, while the conditions $k_\alpha\in\omega_+\omega_-\mathbb{C}_+$ are
needed for the convergence of the integral. The requirement $\omega_+\omega_-\in\mathbb{H}_-$ is imposed to end up with an
arcwise connected parameter space.
\end{rema}


\section{The proof of the integral identity}
 
Since $\mathcal{S}$ is arcwise connected, it suffices to prove Theorem \ref{hyperbolicCM} for parameters in the smaller parameter domain
\[
\mathcal{S}^\prime:=\{(\omega_+,\omega_-,k)\in (\mathbb{C}_+\cap\mathbb{H}_-)^{\times 2}\times\mathcal{K} \,\, | \,\, 
\omega_+/\omega_-\in\mathbb{H}_+,\,\,\,\, k_\alpha\in
\mathbb{C}_-\cap\omega_+\mathbb{C}_+\cap
(\omega_-\mathbb{H}_+-\omega_+)\,\,\,\forall\alpha\in\Sigma\}
\]
(hereby thus destroying the symmetric role of $\omega_+$ and $\omega_-$). Note that for $\omega_{\pm}\in\mathbb{C}_+\cap\mathbb{H}_-$
with $\omega_+/\omega_-\in\mathbb{H}_+$ the allowed space $\mathbb{C}_-\cap\omega_+\mathbb{C}_+\cap(\omega_-\mathbb{H}_+-\omega_+)$ in which
the multiplicity functions may take their values is an open triangle in the third quadrant $\mathbb{C}_-\cap\mathbb{H}_-$ of the complex plane.

We assume the parameter conditions $(\omega_+,\omega_-,k)\in\mathcal{S}^\prime$ throughout this
section. The bases associated to $(\omega_+,u_\alpha\omega_-)$ for $\alpha\in\Sigma$ are denoted by
\[q_\alpha=\exp(2\pi i\omega_+/u_\alpha\omega_-),\qquad \widetilde{q}_\alpha=\exp(-2\pi iu_\alpha\omega_-/\omega_+),
\]
which, under the present assumptions, have moduli $<1$. We also write
\[t_\alpha=\exp(-2\pi ik_\alpha/u_\alpha\omega_-),\qquad \widetilde{t}_\alpha=\exp(-2\pi ik_\alpha/\omega_+),
\]
which also have moduli $<1$. Finally it is convenient to use the notation
\[q_\alpha^z=\exp(2\pi i\omega_+ z/u_\alpha\omega_-),\qquad \widetilde{q}_\alpha^z=\exp(-2\pi iu_\alpha\omega_- z/\omega_+)
\]
for $z\in\mathbb{C}$, 
which allows us to write $t_\alpha=q_\alpha^{-k_\alpha/\omega_+}$ and $\widetilde{t}_\alpha=\widetilde{q}_\alpha{}^{k_\alpha/u_\alpha\omega_-}$. 

\subsection{Splitting the integral}\label{S31}

The starting point for the proof of Theorem \ref{hyperbolicCM} is the following observation.
\begin{lem}\label{Cauchy}
We have 
\[\int_VI(v)dv=(i\omega_-)^n\int_VI(i\omega_- v)dv, 
\]
with both sides absolutely convergent.
\end{lem}
\begin{proof}
Write $\vartheta\in(0,\pi/2)$ for the argument of $i\omega_-\in\mathbb{C}_+\cap\mathbb{H}_+$.
Set $E_{\vartheta}^+=\{re^{i\theta} \, | \, r\geq 0,\,\, 0\leq \theta\leq\vartheta \}$, which is the closure of the open
wedge $\omega_-\mathbb{C}_+\cap\mathbb{H}_+$ in $\mathbb{C}$, and write
\[V_\vartheta^+:=\{v\in V_{\mathbb{C}} \, | \, \langle \alpha,v\rangle\in E_{\vartheta}^+,\,\,\, \forall \alpha\in\Sigma^+ \}.
\]
Observe that $v=\sum_{j=1}^n\lambda_j\omega_j^\vee\in V_\vartheta^+$ iff $\lambda_j\in E_{\vartheta}^+$ for $1\leq j\leq n$.
We apply Cauchy's Theorem to rotate the coordinate-wise integrations over the half lines $[0,\infty)$ to $e^{i\vartheta}[0,\infty)$ in the integral
\[\int_{V^+}I(v)dv=\int_{\lambda_1=0}^\infty\cdots \int_{\lambda_n=0}^{\infty}I\Bigl(\sum_{j=1}^n\lambda_j\omega_j^\vee\Bigr)d\lambda_1\cdots d\lambda_n.
\]
To justify the application of Cauchy's Theorem, we have to show that the integrand $I(v)$ is analytic at $v\in V_{\vartheta}^+$ and that 
$|I(v)|$ has sufficient uniform asymptotic decay when $|v|\rightarrow\infty$ for $v\in V_{\vartheta}^+$.
 
Since the poles of the integrand $I(v)$
(see \eqref{integrand}) are at $\langle\alpha^\prime,v\rangle\in -ik_\alpha+\Lambda_\alpha$ ($\alpha\in\Sigma$) and 
\[\bigl(-ik_\alpha+\Lambda_\alpha\bigr)\cap \bigl(E_{\vartheta}^+\cup (-E_{\vartheta}^+)\bigr)=\emptyset,
\qquad\forall\,\alpha\in\Sigma,
\]
the integrand $I(v)$ is analytic at $v\in V_{\vartheta}^+$. Unfortunately, the asymptotic
estimates for the hyperbolic gamma function $G(z)$ from \cite[Appendix A]{R} and \cite[Cor. 2.3]{Rains} are not good enough to establish the necessary
uniform bounds on the integrand $I(v)$ for $v\in V_{\vartheta}^+$. In fact, for our purposes we would need uniform asymptotics of $G(z)$ for $z$ in suitable translates of $E_{\vartheta}^+$,  
but these regions are not allowed in \cite[Cor. 2.3]{Rains} because $E_{\vartheta}^+$ includes the half-line $\mathbb{R}_{\geq 0}e^{i\vartheta}$ running parallel to the wedge 
$i\omega+\mathbb{R}_{\geq 0}i\omega_++\mathbb{R}_{\geq 0}i\omega_-$ generated by the zeros $i\omega+\Lambda$ of $G(z)$. To bypass this
problem, we establish the necessary asymptotics of $I(v)$ using Shintani's product formula \eqref{ShG} for $G(z)$, which we are allowed to use since
$\omega_+/\omega_-\in\mathbb{H}_+$.

Using the alternative expression \eqref{integrandalternative} for the integrand $I(v)$ in combination with
the reflection equation \eqref{reflectionequation} and the product formula \eqref{ShG} for the hyperbolic gamma function, we then write
\begin{equation}\label{Iasym}
I(v)=\exp\left(-\frac{4\pi\langle\rho_k,v\rangle}{\omega_+\omega_-}\right)\prod_{\alpha\in\Sigma^+}\frac{A_+^\alpha\bigl(\phi_+(\langle\alpha^\prime,v\rangle)\bigr)A_-^\alpha\bigl(
\phi_-(\langle\alpha,v\rangle)\bigr)}
{B_+^\alpha\bigl(\phi_+(\langle\alpha^\prime,v\rangle)\bigr)B_-^\alpha\bigl(\phi_-(\langle\alpha,v\rangle)\bigr)}
\end{equation}
with $\phi_{\pm}(z)=\exp(-2\pi z/\omega_{\pm})$ and with the four complex analytic functions
\begin{equation*}
\begin{split}
A_+^\alpha(z)&=(1-z)\bigl(\widetilde{q}_\alpha\widetilde{t}_\alpha z;\widetilde{q}_\alpha\bigr)_{\infty},\quad
A_-^\alpha(z)=(1-z)\bigl(q_\alpha t_\alpha^{-1}z;q_\alpha\bigr)_{\infty},\\
B_+^\alpha(z)&=\bigl(\widetilde{t}_\alpha^{-1}z;\widetilde{q}_\alpha\bigr)_{\infty},\qquad\qquad\,\,\,\,
B_-^\alpha(z)=\bigl(t_\alpha z;q_\alpha\bigr)_{\infty}.
\end{split}
\end{equation*}
Our goal is to show that the $\alpha$-dependent factors in \eqref{Iasym} are uniformly bounded on $V_{\vartheta}^+$.
Fix $\alpha\in\Sigma^+$, then $V_{\vartheta}^+$ is mapped onto $E_{\vartheta}^+$ by the complex linear functionals $\langle \alpha, \cdot\rangle$ and 
$\langle \alpha^\prime, \cdot\rangle$. By the parameter conditions we have $E_{\vartheta}^+\subset \omega_{\pm}\overline{\mathbb{C}}_+$ (where $\overline{\mathbb{C}}_+$ is the closed right
half plane), hence $\phi_{\pm}$ map $E_{\vartheta}^+$ into the closed unit disc $D=\{z\in\mathbb{C} \,\, | \,\, |z|\leq 1\}$.
Consequently, $A_{\pm}^\alpha\circ \phi_{\pm}, B_{\pm}^\alpha\circ \phi_{\pm}: E_{\vartheta}^+\rightarrow \mathbb{C}$ are bounded. 
Furthermore, $B_-^\alpha$ is zero-free on $D$ since $|t_\alpha|<1$, hence $(B_-^\alpha\circ \phi_-)^{-1}$ is bounded on $E_{\vartheta}^+$. 
Now $|\widetilde{t}_\alpha^{-1}|>1$, so $B_+^\alpha$ has a finite number of poles 
$\{\widetilde{t}_\alpha\widetilde{q}_\alpha^{-j}\}_{j=0}^r$ in $D$, where $r$ is the largest nonnegative integer such that $\widetilde{t}_\alpha\widetilde{q}_\alpha^{-r}\in D$.
Yet we still claim that $E_{\vartheta}^+\ni z\mapsto B_+^\alpha(\phi_+(z))$ is zero free and that its inverse is bounded on $E_{\vartheta}^+$. For this it suffices to show that
the poles $\{\widetilde{t}_\alpha\widetilde{q}_\alpha^{-j}\}_{j=0}^r$ of $B_+^\alpha$ in $D$ are not contained in the closure of $\phi_+(E_{\vartheta}^+)$.

Fix $0\leq j\leq r$ and choose an open neighborhood $U$ of $k_\alpha$ contained in the open triangle $\mathbb{C}_-\cap\omega_+\mathbb{C}_+\cap (\omega_-\mathbb{H}_+-\omega_+)$.
Since $\phi_+$ is an open map, $\widetilde{q}_\alpha^{-j}\phi_+(iU)$ is an open neighborhood of $\widetilde{q}_\alpha^{-j}\widetilde{t}_\alpha$. It suffices to show that 
$\phi_+(E_{\vartheta}^+)\cap\widetilde{q}_\alpha^{-j}\phi_+(iU)=\emptyset$, or equivalently that $E_{\vartheta}^+\cap U_j(l)=\emptyset$ for all $l\in\mathbb{Z}$, where
$U_j(l)=i(U-ju_\alpha\omega_-+l\omega_+)$. 

If $l\leq 0$ then $U_j(l)\subset \mathbb{H}_-$ since $U\subset \mathbb{C}_-$ and $-\omega_{\pm}\in\mathbb{C}_-$. On the other hand, $E_{\vartheta}^+$
is contained in the closed upper half plane, hence $E_{\vartheta}^+\cap U_j(l)=\emptyset$. 

If $l>0$ then $U-ju_\alpha\omega_-+l\omega_+\subset\omega_-\mathbb{H}_+$,
hence $U_j(l)\subset \omega_-\mathbb{C}_-$. On the other hand, $E_{\vartheta}^+$ is contained in the closure of $\omega_-\mathbb{C}_+$, hence again
$E_{\vartheta}^+\cap U_j(l)=\emptyset$, as required.

The bounds thus obtained now give, in combination with \eqref{Iasym}, the uniform estimate
\[|I(v)|\leq C\left|\exp\left(-\frac{4\pi\langle\rho_k,v\rangle}{\omega_+\omega_-}\right)\right|,\qquad \forall\,v\in V_{\varphi}^+
\]
for some constant $C>0$. 
In particular, for $v=\sum_{j=1}^n\lambda_j\omega_j^\vee\in V_{\vartheta}^+$ we have
\[|I(v)|\leq C\prod_{j=1}^n\exp(4\pi c_j|\lambda_j|)
\]
with
\[c_j=\underset{0\leq\theta\leq\vartheta}{\hbox{Max}}\hbox{Re}\left(-\frac{\langle\rho_k,\omega_j^\vee\rangle e^{i\theta}}{\omega_+\omega_-}\right).
\]
By the parameter conditions we have $c_j<0$ for $1\leq j\leq n$, hence Cauchy's theorem can be applied repeatedly to obtain
\[\int_{V^+}I(v)dv=(i\omega_-)^n\int_{V^+}I(i\omega_- v)dv.
\]
Symmetrizing both sides yields the desired result.
\end{proof}
The product formula \eqref{ShG} applied to the expression \eqref{integrand} of the integrand $I(v)$ gives
\[I(i\omega_- v)=K\Delta(v)\widetilde{\Delta}(v)
\]
with the constant 
\begin{equation}\label{K}
K=\prod_{\alpha\in\Sigma^+}\exp\bigl(-\pi ik_\alpha(k_\alpha+2\omega_\alpha)/u_\alpha\omega_+\omega_-\bigr)
\end{equation}
and with
\begin{equation}\label{Delta}
\begin{split}
\Delta(v)&=\prod_{\alpha\in\Sigma}\frac{\bigl(\exp(2\pi i\langle\alpha,v\rangle);q_\alpha\bigr)_{\infty}}
{\bigl(t_\alpha\exp(2\pi i\langle\alpha,v\rangle);q_\alpha\bigr)_{\infty}},\\
\widetilde{\Delta}(v)&=
\prod_{\alpha\in\Sigma}\frac{\bigl(\widetilde{t}_\alpha\widetilde{q}_\alpha^{1+\langle\alpha,v\rangle};\widetilde{q}_\alpha\bigr)_{\infty}}
{\bigl(\widetilde{q}_\alpha^{1+\langle\alpha,v\rangle};\widetilde{q}_\alpha\bigr)_{\infty}}.
\end{split}
\end{equation}
Observe that $\Delta(v)$ is $P^\vee$-invariant. By Lemma \ref{Cauchy} 
we can thus write
\begin{equation}\label{splitintegral}
\int_VI(v)dv=(i\omega_-)^nK\int_D\Delta(v)\left(\sum_{\lambda\in P^\vee}\widetilde{\Delta}(v+\lambda)\right)dv,
\end{equation}
which plays a crucial role in the derivation of the integral identity \eqref{CMalternative} in the following subsection.


\subsection{Cherednik-Macdonald identities}\label{S32}

We prove the desired integral identity
\begin{equation}\label{CMalternative2}
\int_{V}I(v)dv=\#W f\omega_-^n\prod_{\alpha\in\Sigma^+}\frac{G_\alpha(i(\langle\rho_k,\alpha^\vee\rangle+\omega_\alpha))G_\alpha(i(\langle\rho_k,\alpha^\vee\rangle-\omega_\alpha))}
{G_{\alpha}(i(\langle\rho_k,\alpha^\vee\rangle+k_\alpha+\omega_\alpha))
G_\alpha(i(\langle\rho_k,\alpha^\vee\rangle-k_\alpha+u_\alpha\omega_-\delta_\alpha-\omega_\alpha))}
\end{equation}
using the alternative expression \eqref{splitintegral} for the left hand side.
The right hand side of \eqref{splitintegral} simplifies by the following identity, which is equivalent to Macdonald's summation identity
\cite[\S 7]{Mac2}.

\begin{thm}\label{Macdonaldsum}\cite{Mac2}
For $v\in V_{\mathbb{C}}$ such that $\langle \alpha,v\rangle\not\in\mathbb{Z}$ for all $\alpha\in\Sigma$,
\[\widetilde{N}:=\sum_{\lambda\in P^\vee}\widetilde{\Delta}(v+\lambda)\]
is absolutely convergent and independent of $v$. Moreover, 
\[\widetilde{N}=f\prod_{\alpha\in\Sigma^+}\frac{\bigl(\widetilde{q}_\alpha\widetilde{t}_\alpha\exp(-2\pi i\langle\rho_k,\alpha^\vee\rangle/\omega_+);\widetilde{q}_\alpha\bigr)_{\infty}
\bigl(\widetilde{q}_\alpha^{\delta_\alpha}\widetilde{t}_\alpha^{-1}\exp(-2\pi i\langle\rho_k,\alpha^\vee\rangle/\omega_+);\widetilde{q}_\alpha\bigr)_{\infty}}
{\bigl(\widetilde{q}_\alpha\exp(-2\pi i\langle\rho_k,\alpha^\vee\rangle/\omega_+);\widetilde{q}_\alpha\bigr)_{\infty}
\bigl(\exp(-2\pi i\langle\rho_k,\alpha^\vee\rangle/\omega_+);\widetilde{q}_\alpha\bigr)_{\infty}}.
\]
\end{thm}
\begin{proof}
We give the precise correspondence with \cite[(7.4)]{Mac2}. In our notations, after analytic continuation in the parameters, it reads as follows.
For $\omega_+/\omega_-\in\mathbb{H}_+$, $\kappa\in\mathcal{K}$ with $\kappa_\alpha\in \omega_-\omega_+^{-1}\mathbb{H}_+$ and for
$w\in V_{\mathbb{C}}$ such that $\langle\alpha^\prime,w\rangle\not\in\mathbb{Z}$ for all $\alpha\in\Sigma$, we have
\begin{equation}\label{Macdonaldversion}
\sum_{\lambda\in\Lambda}\prod_{\alpha\in\Sigma}\frac{\bigl(q_\alpha^{1+\kappa_\alpha+\langle\alpha^\prime, w+\lambda\rangle};q_\alpha\bigr)_{\infty}}
{\bigl(q_\alpha^{1+\langle\alpha^\prime, w+\lambda\rangle};q_\alpha\bigr)_{\infty}}
=f\prod_{\alpha\in\Sigma^+}\frac{\bigl(q_\alpha^{\langle\rho_{\kappa},\alpha^\vee\rangle+\kappa_\alpha+1}, q_\alpha^{\langle\rho_\kappa,\alpha^\vee\rangle-\kappa_\alpha+\delta_\alpha};q_\alpha\bigr)_{\infty}}
{\bigl(q_\alpha^{\langle\rho_\kappa,\alpha^\vee\rangle+1}, q_\alpha^{\langle\rho_\kappa,\alpha^\vee\rangle};q_\alpha\bigr)_{\infty}}
\end{equation}
with the sum being absolutely convergent and with $\Lambda=P^\vee$ in case {\bf (i)} (in which case $u_\alpha=1$ and $\alpha^\prime=\alpha$ for all $\alpha\in\Sigma$)
and with $\Lambda=P$ in case {\bf (ii)} (in which case $u_\alpha=2/\|\alpha\|^2$ and $\alpha^\prime=\alpha^\vee$ for all $\alpha\in\Sigma$).

We relate this formula to the statement in the theorem for case {\bf (i)} and case {\bf (ii)} separately.
The theorem for case {\bf (i)} with parameters $(\omega_+^\prime,\omega_-^\prime,k)\in\mathcal{S}^\prime$ 
follows from case {\bf (i)} of \eqref{Macdonaldversion} with $w=v$ and the parameters specialized to $(\omega_+,\omega_-,\kappa)=
(-\omega_-^\prime,\omega_+^\prime,k/\omega_-^\prime)$, where $k/\omega_-^\prime$ is the multiplicity function that takes value
$k_\alpha/\omega_-^\prime$ at $\alpha\in\Sigma$.

For case {\bf (ii)} we write $\varphi\in\Sigma$ for the highest root with respect to $\Sigma^+$ (which is a long root).
For $\alpha\in\Sigma$ we define $\widetilde{\alpha}=u_{\varphi}^{-\frac{1}{2}}\alpha^\vee$. Then $\widetilde{\Sigma}=\{\widetilde{\alpha}\}_{\alpha\in\Sigma}$
is a reduced irreducible root system in $V$, normalized so that short roots have squared length two. We take $\{\widetilde{\alpha}_i\}_{i=1}^n$ as basis
of $\widetilde{\Sigma}$. Observe that the weight lattice of $\widetilde{\Sigma}$ is $u_{\varphi}^{-\frac{1}{2}}P^\vee$ and that $u_{\widetilde{\alpha}}=u_\varphi/u_\alpha$.
The theorem for case {\bf (ii)} with parameters $(\omega_+^\prime,\omega_-^\prime,k)\in\mathcal{S}^\prime$ now follows from case {\bf (ii)} of \eqref{Macdonaldversion} with $\Sigma$ replaced by
$\widetilde{\Sigma}$, $w=u_\varphi^{-\frac{1}{2}}v$, $(\omega_+,\omega_-)=(-u_\varphi\omega_-^\prime,\omega_+^\prime)$ and $\kappa_{\widetilde{\alpha}}=k_\alpha/u_\alpha\omega_-^\prime$.
\end{proof}

Theorem \ref{Macdonaldsum} and \eqref{splitintegral} thus yield
\begin{equation}\label{onefactor}
\int_VI(v)dv=(i\omega_-)^nK\widetilde{N}\int_D\Delta(v)dv.
\end{equation}
The following constant term identity, which was conjectured by Macdonald \cite[\S 12]{Mac1} and proved by
Cherednik in \cite{C} using double affine Hecke algebras, gives the evaluation of the resulting integral.
\begin{thm}\label{Macdonaldintegral}\cite{C, Mac1}
We have
\[\int_D\Delta(v)dv=N
\]
with
\[N=\#W\prod_{\alpha\in\Sigma^+}\frac{\bigl(\exp(-2\pi i\langle\rho_k,\alpha^\vee\rangle/u_\alpha\omega_-);q_\alpha\bigr)_{\infty}
\bigl(q_\alpha\exp(-2\pi i\langle\rho_k,\alpha^\vee\rangle/u_\alpha\omega_-);q_\alpha\bigr)_{\infty}}
{\bigl(t_\alpha\exp(-2\pi i\langle\rho_k,\alpha^\vee\rangle/u_\alpha\omega_-);q_\alpha\bigr)_{\infty}\bigl(q_\alpha t_\alpha^{-1}\exp(-2\pi i\langle\rho_k,\alpha^\vee\rangle/u_\alpha\omega_-);q_\alpha\bigr)_{\infty}}.
\]
\end{thm}
Combining Theorem \ref{Macdonaldintegral} and \eqref{onefactor} we thus obtain the integral evaluation
\begin{equation}\label{zerofactor}
\int_VI(v)dv=(i\omega_-)^nKN\widetilde{N}.
\end{equation}
By the product formula \eqref{ShG} for the hyperbolic gamma function we have
\[(i\omega_-)^nKN\widetilde{N}=\#W f\omega_-^n\prod_{\alpha\in\Sigma^+}\frac{G_\alpha(i(\langle\rho_k,\alpha^\vee\rangle+\omega_\alpha))G_\alpha(i(\langle\rho_k,\alpha^\vee\rangle-\omega_\alpha))}
{G_\alpha(i(\langle\rho_k,\alpha^\vee\rangle+k_\alpha+\omega_\alpha))G_\alpha(i(\langle\rho_k,\alpha^\vee\rangle-k_\alpha+u_\alpha\omega_-\delta_\alpha-\omega_\alpha))},
\]
hence \eqref{zerofactor} gives the desired integral identity \eqref{CMalternative2}. This completes the proof of 
Theorem \ref{hyperbolicCM}.


\section{Connection to type II integral evaluations}\label{s4}

In this section we show that the generalized Cherednik-Macdonald constant term identity for root system $\Sigma$ of type
$A_1$, $B_n$ and $C_n$ is a special case of the first level degeneration \cite[Thm. 5]{vDS} of the type II multivariate hyperbolic
integral evaluation \cite[Thm. 4]{vDS}, see also \cite[Cor. 4.4]{Rains}. This integral identity is the hyperbolic analogue of
the constant term associated to Koornwinder's \cite{K} extension of the Macdonald theory. Before giving the identity we first have to 
introduce some more notations.

Consider $\mathbb{R}^n$  with standard scalar product $\langle\cdot,\cdot\rangle$ and standard orthonormal basis $\{\epsilon_j\}_{j=1}^n$.
We 
write $d^\prime v$ for the Lebesgue measure on $\mathbb{R}^n$ 
normalized by $\int_{[0,1]^n}d^\prime v=1$.
Set 
\[R_s=\{\pm \epsilon_j\}_{j=1}^n,\qquad R_l=\{\pm(\epsilon_r\pm\epsilon_s)\}_{1\leq r<s\leq n},
\]
where $R_l$ is the empty-set if $n=1$ (all sign combinations are allowed).
We write $\underline{\gamma}=(\gamma_1,\gamma_2,\gamma_3,\gamma_4)\in\mathbb{C}^{4}$ and $|\underline{\gamma}|=\gamma_1+\gamma_2+\gamma_3+\gamma_4$. 
Define the parameter space
$\mathcal{S}_{BC}$ to be the set of parameters $(\omega_+,\omega_-,\underline{\gamma},\kappa)\in\mathbb{C}_+^{2}\times\mathbb{C}^{5}$
satisfying 
\begin{equation}\label{conditions}
\begin{split}
\omega_+\omega_-&\in\mathbb{H}_-,\\
\gamma_r&\in\mathbb{C}_-\qquad\qquad (1\leq r\leq 4),\\
\omega_++\omega_-+|\underline{\gamma}|&\in\omega_+\omega_-\mathbb{C}_+,\\
\kappa&\in\mathbb{C}_-\cap\omega_+\omega_-\mathbb{C}_+. 
\end{split}
\end{equation}
It serves as the analogue of the parameter space $\mathcal{S}$ for the generalized Cherednik-Macdonald identities.
For $n=1$ the space $\mathcal{S}_{BC}$ should be interpreted as the set of parameters $(\omega_+,\omega_-,\underline{\gamma})\in\mathbb{C}_+^{2}\times\mathbb{C}^{4}$
that are subject to the conditions given by the first three lines of \eqref{conditions}.

For $(\omega_+,\omega_-,\underline{\gamma},\kappa)\in\mathcal{S}_{BC}$ we now consider the integral
\[
J_{BC}:=\int_{\mathbb{R}^n}\prod_{\alpha\in R_s}\frac{G(\langle\alpha,v\rangle+i\omega)G(\langle\alpha,v\rangle+\frac{i\omega_+}{2})
G(\langle\alpha,v\rangle+\frac{i\omega_-}{2})}{\prod_{j=1}^4G(\langle\alpha,v\rangle+i(\omega+\gamma_j))}
\prod_{\beta\in R_l}\frac{G(\langle\beta,v\rangle+i\omega)}{G(\langle\beta,v\rangle+i(\omega+\kappa))}d^\prime v.
\]
Note that the integrand can be simplified since
\[\prod_{\alpha\in R_s}G(2\langle\alpha,v\rangle+i\omega)=
\prod_{\alpha\in R_s}G(\langle\alpha,v\rangle+i\omega)G\bigl(\langle\alpha,v\rangle+\frac{i\omega_+}{2}\bigr)
G\bigl(\langle\alpha,v\rangle+\frac{i\omega_-}{2}\bigr),
\]
which follows from the identity
\begin{equation}\label{automorphy}
G(2z+i\omega)=G(z)G\bigl(z+\frac{i\omega_+}{2}\bigr)G\bigl(z+\frac{i\omega_-}{2}\bigr)G(z+i\omega)
\end{equation}
(see \cite[Prop. III.2]{R0}) and the reflection equation \eqref{reflectionequation}. The integral $J_{BC}$ is then seen to be equivalent
to the integral from \cite[Thm 5]{vDS}. It follows from \cite[Thm. 5]{vDS} that $J_{BC}$ is 
absolutely convergent for $(\omega_+,\omega_-,\underline{\gamma},\kappa)\in\mathcal{S}_{BC}$ and that it is equal to 
\begin{equation}\label{NBC}
N_{BC}=2^nn!\bigl(\sqrt{\omega_+\omega_-}\bigr)^n\prod_{j=0}^{n-1}\frac{G(i(\omega+\kappa))G(i(\omega+(2n-j-2)\kappa+|\underline{\gamma}|))}
{G(i(\omega+(j+1)\kappa))\prod_{1\leq r<s\leq 4}G(i(\omega+j\kappa+\gamma_r+\gamma_s))}.
\end{equation}
This can also be proved by a straightforward generalization of the arguments in this paper, replacing the role of the Cherednik-Macdonald
constant term identity by Gustafson's \cite{Gus} multivariate
Askey-Wilson integral evaluation and the role of Macdonald's summation identity by van Diejen's \cite{vDi} multivariate ${}_6\Psi_6$ summation formula. 
It can also be obtained as rigorous limit from the multivariate hyperbolic integral evaluation of type II (\cite[Thm. 4]{vDS}, \cite{Rains}),
see \cite{BR}. For $n=1$, the integral identity is the hyperbolic Askey-Wilson integral evaluation from \cite{R1, S}. 

We now proceed to show that the generalized Cherednik-Macdonald integrals $\int_VI(v)dv$ for root systems $\Sigma$ of types $A_1$, $B_n$
and $C_n$ are special cases of the hyperbolic integral $J_{BC}$. 
We leave the (rather cumbersome) identification of their exact evaluations
to the evaluation formula $N_{BC}$ of $J_{BC}$ to the reader. For this identification one needs the identities $G(0)=1$ (which is obvious) and
\[
G\Bigl(\frac{i\omega_{\pm}}{2}\Bigr)=\exp\Bigl(\frac{1}{2}\int_0^{\infty}dy\Bigl(\frac{1}{\omega_{\mp} y^2}-\frac{1}{y\sinh(\omega_{\mp}y)}\Bigr)\Bigr)=
\sqrt{2},
\]
see e.g. \cite[(3.17)]{R0} for the last equality.

\subsection{$\Sigma$ of type $A_1$ and of type $B_n$}
We realize the root system $\Sigma$ of type $B_n$ ($n\geq 2$) as $\Sigma=\Sigma_s\cup\Sigma_l\subset\mathbb{R}^n$ with
\[\Sigma_s=\{\pm\sqrt{2}\epsilon_j\}_{j=1}^n,\qquad \Sigma_l=\{\pm\sqrt{2}(\epsilon_r\pm\epsilon_s)\}_{1\leq r<s\leq n},
\]
so $\Sigma_s=\sqrt{2}R_s$ and $\Sigma_l=\sqrt{2}R_l$
(the awkward normalization is caused by the requirement that the short roots $\Sigma_s$ in $\Sigma$ should have squared length two).
We include $n=1$ as the case $\Sigma=\Sigma_s=\{\pm\sqrt{2}\epsilon_1\}$, in which case it is the root system of type $A_1$.
In the discussion below the obvious modifications have to be made for $n=1$.
We have $P^\vee=\frac{1}{\sqrt{2}}\bigoplus_{j=1}^n\mathbb{Z}\epsilon_j$, so 
$dv=(\sqrt{2})^nd^\prime v$. 
For a given multiplicity function $k\in\mathcal{K}$
we write $k_s$ for its value on the short roots $\Sigma_s$ and $k_l$ for its value on the long roots $\Sigma_l$.\\

Case {\bf (i)} ($u_\alpha=1$ for all $\alpha\in\Sigma$): After a change of integration variables the generalized Cherednik-Macdonald
integral becomes
\[\int_VI(v)dv=\int_{\mathbb{R}^n}\prod_{\alpha\in R_s}\frac{G(\langle\alpha,v\rangle+i\omega)}{G(\langle \alpha,v\rangle+i(k_s+\omega))}
\prod_{\beta\in R_l}\frac{G(\langle\beta,v\rangle+i\omega)}{G(\langle\beta,v\rangle+i(k_l+\omega))}d^\prime v
\]
for $(\omega_+,\omega_-,k)\in\mathcal{S}$. It equals $J_{BC}$ with the parameters $(\underline{\gamma},\kappa)$
specialized to
\[(\underline{\gamma},\kappa)=\bigl(k_s,-\frac{\omega_+}{2}, -\frac{\omega_-}{2}, -\omega, k_l\bigr).
\]
Under this parameter correspondence the requirement $(\omega_+,\omega_-,\underline{\gamma},\kappa)\in \mathcal{S}_{BC}$ exactly corresponds to
$(\omega_+,\omega_-,k)\in\mathcal{S}$. 

In particular, with this parameter specialization, the hyperbolic Askey-Wilson integral evaluation reduces to the generalized Cherednik-Macdonald identity for $\Sigma$ of type $A_1$:
\begin{equation}\label{A1}
\int_{\mathbb{R}}\frac{G(v+i\omega)G(-v+i\omega)}{G(v+i(k+\omega))G(-v+i(k+\omega))}d^\prime v=
4\sqrt{\omega_+\omega_-}\frac{G(i(k+\omega))G(i(k-\omega))}{G(i(2k+\omega))}
\end{equation}
for $(\omega_+,\omega_-,k)\in\mathcal{S}$.\\

Case {\bf (ii)} ($u_\alpha=2/\|\alpha\|^2$ for $\alpha\in\Sigma$): The generalized Cherednik-Macdonald
integral becomes
\[\int_VI(v)dv=\int_{\mathbb{R}^n}\prod_{\alpha\in R_s}\frac{G(\omega_+,\omega_-;\langle\alpha,v\rangle+i\omega)}{G(\omega_+,\omega_-;\langle \alpha,v\rangle+i(k_s+\omega))}
\prod_{\beta\in R_l}\frac{G\bigl(\omega_+,\frac{\omega_-}{2};\frac{\langle\beta,v\rangle}{2}+i(\frac{\omega_+}{2}+\frac{\omega_-}{4})\bigr)}{G\bigl(\omega_+,\frac{\omega_-}{2};
\frac{\langle\beta,v\rangle}{2}+i(k_l+\frac{\omega_+}{2}+\frac{\omega_-}{4})\bigr)}d^\prime v
\]
for $(\omega_+,\omega_-,k)\in\mathcal{S}$. We rewrite the integrand involving only hyperbolic gamma functions with quasi-periods $(2\omega_+,\omega_-)$ 
using  
\[
G(\omega_+,\omega_-;z)=G\bigl(2\omega_+,\omega_-;z+\frac{i\omega_+}{2}\bigr)G\bigl(2\omega_+,\omega_-;z-\frac{i\omega_+}{2}\bigr)
\]
and $G(\omega_+,\omega_-/2;z)=G(2\omega_+,\omega_-;2z)$ (see \cite[Prop. III.2]{R0}).
We conclude that $\int_VI(v)dv$ equals $J_{BC}$ with respect to the quasi-periods $(2\omega_+,\omega_-)$ and with the parameters $(\underline{\gamma},\kappa)$
specialized to
\[(\underline{\gamma},\kappa)=\bigl(k_s,k_s-\omega_+, -\frac{\omega_-}{2}, -\omega_+-\frac{\omega_-}{2}, 2k_l\bigr).
\]
The requirement $(\omega_+,\omega_-,\underline{\gamma},\kappa)\in \mathcal{S}_{BC}$ then becomes $(\omega_+,\omega_-,k)\in\mathcal{S}$.

\subsection{$\Sigma$ of type $C_n$}

We realize the root system $\Sigma$ of type $C_n$ ($n\geq 2$) as $\Sigma=\Sigma_s\cup\Sigma_l\subset\mathbb{R}^n$ with
\[\Sigma_s=\{\pm (\epsilon_r\pm\epsilon_s)\}_{1\leq r<s\leq n},\qquad \Sigma_l=\{\pm 2\epsilon_j\}_{j=1}^n,
\]
so $\Sigma_s=R_l$ and $\Sigma_l=2R_s$.
We have $f=\#(P^\vee/Q^\vee)=2$ and $Q^\vee=\bigoplus_{j=1}^n\mathbb{Z}\epsilon_j$, so $d^\prime v=\frac{1}{2}dv$.
For a given multiplicity function $k\in\mathcal{K}$
we write $k_s$ for its value on the short roots $\Sigma_s$ and $k_l$ for its value on the long roots $\Sigma_l$.\\

Case {\bf (i)} ($u_\alpha=1$ for all $\alpha\in\Sigma$):
The generalized Cherednik-Macdonald integral becomes
\[
\int_VI(v)dv=2\int_{\mathbb{R}^n}\prod_{\alpha\in R_s}\frac{G(2\langle\alpha,v\rangle+i\omega)}
{G\bigl(2(\langle\alpha,v\rangle+\frac{ik_l}{2})+i\omega\bigr)}\prod_{\beta\in R_l}\frac{G(\langle\beta,v\rangle+i\omega)}{G(\langle\beta,v\rangle+i(k_s+\omega))}
d^\prime v
\]
for $(\omega_+,\omega_-,k)\in\mathcal{S}$. By \eqref{automorphy} and \eqref{reflectionequation}, the product over $\alpha\in R_s$ can be rewritten as
\[
\prod_{\alpha\in R_s}\frac{G(\langle\alpha,v\rangle+i\omega)G\bigl(\langle\alpha,v\rangle+\frac{i\omega_+}{2}\bigr)G\bigl(\langle\alpha,v\rangle+\frac{i\omega_-}{2}\bigr)}
{G\bigl(\langle \alpha,v\rangle+i(\frac{k_l}{2}+\omega)\bigr)G\bigl(\langle \alpha,v\rangle+i(\frac{k_l}{2}+\frac{\omega_+}{2})\bigr)
G\bigl(\langle \alpha,v\rangle+i(\frac{k_l}{2}+\frac{\omega_-}{2})\bigr)G\bigl(\langle \alpha,v\rangle+\frac{ik_l}{2}\bigr)}.
\]
Hence $\int_VI(v)dv$ equals $2J_{BC}$ with the parameters $(\underline{\gamma},\kappa)$
specialized to
\[(\underline{\gamma},\kappa)=\Bigl(\frac{k_l}{2}, \frac{k_l}{2}-\frac{\omega_+}{2}, \frac{k_l}{2}-\frac{\omega_-}{2}, \frac{k_l}{2}-\omega, k_s\Bigr).
\]
The parameter conditions $(\omega_+,\omega_-,\underline{\gamma},\kappa)\in \mathcal{S}_{BC}$ become $(\omega_+,\omega_-,k)\in\mathcal{S}$.

Case {\bf (ii)} ($u_\alpha=2/\|\alpha\|^2$ for all $\alpha\in\Sigma$): The generalized Cherednik-Macdonald
integral becomes
\[
\int_VI(v)dv=2\int_{\mathbb{R}^n}\prod_{\alpha\in R_s}\frac{G\bigl(\omega_+,\frac{\omega_-}{2}; \langle\alpha,v\rangle+i(\frac{\omega_+}{2}+\frac{\omega_-}{4})\bigr)}
{G\bigl(\omega_+,\frac{\omega_-}{2};\langle\alpha,v\rangle+i(k_l+\frac{\omega_+}{2}+\frac{\omega_-}{4})\bigr)}
\prod_{\beta\in R_l}\frac{G(\omega_+,\omega_-;\langle\beta,v\rangle+i\omega)}{G(\omega_+,\omega_-;\langle\beta,v\rangle+i(k_s+\omega))}d^\prime v.
\]
We rewrite the integrand involving only hyperbolic gamma functions with quasi-periods $(\omega_+,\omega_-)$ 
using  
\[
G\bigl(\omega_+,\frac{\omega_-}{2};z\bigr)=G\bigl(\omega_+,\omega_-;z+\frac{i\omega_-}{4}\bigr)G\bigl(\omega_+,\omega_-;z-\frac{i\omega_-}{4}\bigr),
\]
cf. \cite[Prop. III.2]{R0}.
We conclude that $\int_VI(v)dv$ equals $2J_{BC}$ with the parameters $(\underline{\gamma},\kappa)$
specialized to
\[(\underline{\gamma},\kappa)=\Bigl(k_l, k_l-\frac{\omega_-}{2}, -\frac{\omega_+}{2}, -\omega, k_s\Bigr).
\]
The parameter conditions $(\omega_+,\omega_-,\underline{\gamma},\kappa)\in \mathcal{S}_{BC}$ again become $(\omega_+,\omega_-,k)\in\mathcal{S}$.


\end{document}